\documentclass[a4paper,twoside]{article}

\usepackage{algorithm,algpseudocode,caption}
\usepackage{array}
\usepackage{dsfont}
\usepackage{multirow}

\usepackage{epsfig}
\usepackage{subcaption}
\usepackage{calc}
\usepackage{amssymb}
\usepackage{amstext}
\usepackage{amsmath}
\usepackage{amsthm}
\usepackage{hyperref}
\usepackage{multicol}
\usepackage{nicefrac}
\usepackage{pslatex}
\usepackage{apalike}
\usepackage{SCITEPRESS}     

\renewcommand{\v}{\vspace{5mm}}

\newcommand{\ie}{\textit{id est} }

\newcommand{\dt}{\delta t}

\newcommand{\dsp}{\displaystyle}
\newcommand{\ndsp}{\textstyle}

\newcommand{\tinit}{t^{[\text{init}]}}
\newcommand{\tend}{t^{[\text{end}]}}

\newcommand{\atol}{\ensuremath{tol_{abs}}}
\newcommand{\rtol}{\ensuremath{tol_{rel}}}

\newcommand{\nsys}{n_{sys}}
\newcommand{\nst}{n_{st}}

\newcommand{\nstk}[1]{n_{st,#1}}
\newcommand{\nink}[1]{n_{in,#1}}
\newcommand{\noutk}[1]{n_{out,#1}}

\newcommand{\rnst}{\mathbb{R}^{\nst}}

\newcommand{\rnstk}[1]{\mathbb{R}^{\nstk{k}}}
\newcommand{\rnink}[1]{\mathbb{R}^{\nink{k}}}
\newcommand{\rnoutk}[1]{\mathbb{R}^{\noutk{k}}}

\newcommand{\Insys}{[\![1, \nsys]\!]}
\newcommand{\Inink}[1]{[\![1, \nink{#1}]\!]}
\newcommand{\Inoutk}[1]{[\![1, \noutk{#1}]\!]}

\newcommand{\bs}{\begin{scriptsize}}
\newcommand{\es}{\end{scriptsize}}

\newcommand{\ddd}[2]{\dsp{\frac{d#1}{d#2}}}
\newcommand{\dspfrac}[2]{\dsp{\frac{#1}{#2}}}

\newcommand{\dspprod}{\dsp{\prod}}
\newcommand{\ndspfrac}[2]{\ndsp{\frac{#1}{#2}}}
\newcommand{\up}[1]{\textsuperscript{#1}}

\newcommand{\restrict}[2]{{#1}_{\big|_{#2}}}
\newcommand{\fornits}{F$ _3$ORNITS }

\begin{document}

\title{\fornits: a flexible variable step size non-iterative co-simulation method handling  subsytems with hybrid advanced capabilities}

\author{\authorname{Yohan \'EGUILLON\sup{1}\orcidAuthor{0000-0002-9386-4646}, Bruno LACABANNE\sup{1}\orcidAuthor{0000-0003-1790-3663} and Damien TROMEUR-DERVOUT\sup{2}\orcidAuthor{0000-0002-0118-8100}}
\affiliation{\sup{1}Siemens Industry Software, Roanne, France}
\affiliation{\sup{2}Institut Camille Jordan, Universit\'e de Lyon ,UMR5208 CNRS-U.Lyon1, Villeurbanne, France}
\email{\{yohan.eguillon, bruno.lacabanne\}@siemens.com, damien.tromeur-dervout@univ-lyon1.fr}
\vspace{-2mm}
}


\abstract{This paper introduces the \fornits non-iterative co-simulation algorithm in which F$ _3$ stands for the $3$ flexible aspects of the method: flexible polynomial order representation of coupling variables, flexible time-stepper applying variable co-simulation step size rules on subsystems allowing it and flexible scheduler orchestrating the meeting times among the subsystems and capable of asynchronousness when subsystems' constraints requires it. The motivation of the \fornits method is to accept any kind of co-simulation model, including any kind of subsystem, regardless on their available capabilities. Indeed, one the major problems in industry is that the subsystems usually have constraints or lack of advanced capabilities making it impossible to implement most of the advanced co-simulation algorithms on them. The method makes it possible to preserve the dynamics of the coupling constraints when necessary as well as to avoid breaking $C^1$ smoothness at communication times, and also to adapt the co-simulation step size in a way that is robust both to zero-crossing variables (contrary to classical relative error-based criteria) and to jumps. Two test cases are presented to illustrate the robustness of the \fornits method as well as its higher accuracy than the non-iterative Jacobi coupling algorithm (the most commonly used method in industry) for a smaller number of co-simulation steps.}

\onecolumn \maketitle \normalsize \setcounter{footnote}{0} \vfill

\section{\uppercase{Introduction}}
\label{sec:introduction}

\noindent Co-simulation consists in processing a simulation of a modular model, that is to say a model composed of several dynamical subsystems connected together. This is usually motivated by the need to simulate a model with multiphysical parts. Designing a subsystem representing the physics of a given field (electricity, mechanics, fluids, thermodynamic, etc) allows the use of a specific and adapted solver for this field, or even the modelling with a specific third party software. Regarding industrial applications, a modular model is prefered because subsystem providers can focus on a part of the global system without taking the rest into account. Nonetheless, gathering different subsystems is not straightforward: a simulation of the equations of the global system cannot be retrieved. Co-simulation is the field investigating the ways to solve such systems, based on regular data communications between the subsystems solved separately.

The co-simulation method (or co-simulation algorithm) is the rule used to process the simulation on such modular systems. It namely deals with: the determination of the times of the data communications, the way the inputs each subsystem should use at each step are computed, and the way the outputs are used and so on. Many co-simulation algorithms have been established until now \cite{Kubler2000} \cite{Arnold2001} \cite{Gu2004} \cite{Bartel2013} \cite{Sicklinger2014}  \cite{Busch2016} and studied \cite{Li2014}  \cite{Schweizer2016}, with various complexity of implementation, subsystems capabilities requirements, or physical field-specific principles (see also the recent state of art on co-simulation of \cite{Gomes2018survey}). It may appear in some co-simulation algorithms that a time interval should be integrated more than once on one or more subsystems. This is called iterative co-simulation. The systems that need to do so must have the capability to replay a simulation over this time-interval with different data (time to reach, input values, ...). Depending on the model provider, this capability, so-called "rollback", is not always available on every subsystem. When it is available, an iterative algorithm such as \cite{Eguillon2019} can be used. Nevertheless, if at least one subsystem cannot rollback, a non-iterative method has to be chosen. The rollback capability is not available in most industrial models. As the \fornits method presented in this paper is non-iterative, it can be applied to configurations with rollback-less subsystems.

\begin{figure}[!ht]
  \vspace{-0.2cm}
  \centering
   {\epsfig{file = 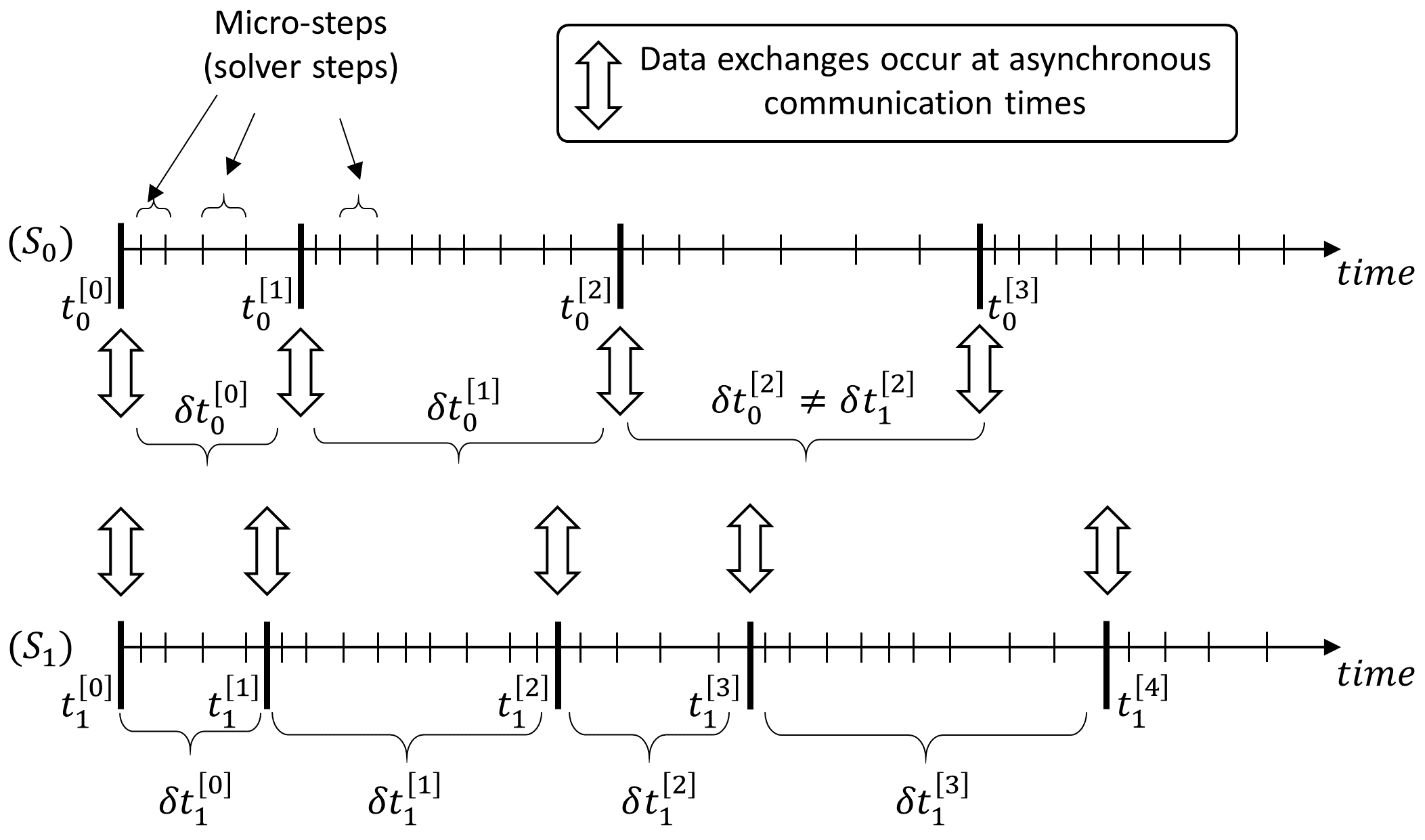, width = 7.5cm}}
  \caption{Visualization of the behavior of a non-iterative asynchronous co-simulation method on $2$ subsystems}
  \label{fig:asynchronous_method}
  \vspace{-0.1cm}
\end{figure}

Besides the rollback, other plateform-dependent capabilities may lead to an impossibility of use of a given co-simulation method on certain subsystems. Amongst them, the ability to provide inputs with $n$\up{th} order time-derivatives and the ability to obtain the time-derivatives of the outputs at the end of a macro-step can be mentioned. The motivation of the \fornits method is to accept any kind of modular model, including any kind of subsystem, regardless on their available capabilities. The first consequence of this specification is that \fornits is non-iterative (in order to accept rollback-less subsystems) and asynchronous (in order to accept modular models with subsystems with an imposed step, even when several subsystems with imposed steps non multiple of one another). Such a method is presented on figure \ref{fig:asynchronous_method}. \fornits also represents the inputs as time-dependent polynomials (when supported by the subsystems concerned), and adapts the communication step size of the subsystems which can handle variable communication step sizes. This both allows accuracy win when frequent exchanges are needed, and saves time when coupling variables can be represented with a high enough accuracy.

This algorithm is based on a variable order polynomial representation of every input (the polynomial degree might be different for each input variable at a same time) determined by an \textit{a posteriori} criterion and redefined at each communication step. A smoother version can be triggered by interpolating on the extrapolated values, as done in \cite{Busch2019} with the so-called "EXTRIPOL" technique. This technique avoids the non-physical jumps on coupling variables at each communication time, and may help the solvers of the subsystems to restart faster after a discontinuity since we can guarantee the $C^1$ smoothness of the input variables. The \fornits method adapts this smoothening with flexible order polynomials, variable step size subsystems, and asynchronous cases.

The paper is structured as follows. Section \ref{sec:mathematical_formalism_and_motivations} presents the  mathematical formalism we used to develop the co-simulation algorithm. This formalism is made of general notations, of subsystems topologies detailing and of common polynomial techniques used in \fornits method. Section \ref{sec:fornits_algorithm} presents the \fornits algorithm. On the one side, it shows the way time-dependent inputs are determined. On the other side the way the step size determination is handled. Section \ref{sec:results_and_behavior_on_two_test_cases} gives the results of \fornits method both on a controlled speed model and on the classical linear two-mass oscillator test case \cite{Schweizer2016} \cite{Eguillon2019}. Comparison of \fornits algorithm with different options and basic non-iterative Jacobi method (most basic non-iterative co-simulation technique with fixed step size) is also achieved. The conclusion is given in section \ref{sec:conclusions}.

\section{\uppercase{Mathematical formalism and motivations}}
\label{sec:mathematical_formalism_and_motivations}

\noindent A subsystem that communicates with other subsystems to form a co-simulation configuration will be represented by its equation. As the context is about time integration, these equations will be time differential equations. First of all, we present the general form of a monolithic system, in other words a closed system (which neither has inputs nor outputs) represented by a $1$\up{st} order differential equation which covers, among others things, ODE, DAE and IDE as follows:
\vspace{-0.2cm}

\begin{equation}
\label{eq:general_diff_equa}
\left\{
	\begin{array}{lll}
		F(\frac{d}{dt}x, x, t) & = & 0\\
		x(\tinit) & = & x^{[\text{init}]}\\
	\end{array}
\right.
\end{equation}

\noindent where
\vspace{-0.5cm}

\begin{equation}
\label{eq:general_diff_equa_decl}
\begin{array}{l}
	\nst \in \mathbb{N}^*, \
	x^{[\text{init}]} \in \mathbb{R}^{\nst}, \
	t\in[\tinit, \tend], \\
	{[}\tinit, \tend] \subset \mathbb{R} \text{ so that } \tend-\tinit \in \mathbb{R}_+^*, \\
	F:\mathbb{R}^{\nst} \times \mathbb{R}^{\nst} \times [\tinit, \tend] \rightarrow \mathbb{R}^{\nst}, \\
\end{array}
\end{equation}

\noindent are given, and
\vspace{-0.3cm}

\begin{equation}
\label{eq:general_diff_equa_decl_time}
x:[\tinit, \tend]\rightarrow\mathbb{R}^{\nst}
\end{equation}
\vspace{-0.2cm}

\noindent is the solution state vector whose $\nst$ components are called  the \textbf{state variables}.

This paper will only cover the ODEs case, which are of the form \eqref{eq:ODE}. The terms "equation", "equations", and "differential equations" will refer to the ODE of a system throughout this document.
\vspace{-0.2cm}

\begin{equation}
\label{eq:ODE}
\frac{d}{dt}x = f(t, x)\\
\end{equation}

\vspace{-0.3cm}
\noindent where
\vspace{-0.3cm}

\begin{equation}
\label{eq:ODE_decl}
f:[\tinit, \tend]\times\rnst\rightarrow\rnst
\end{equation}

\subsection{Framework and notations}
\label{subsection:framework_and_notations}

\noindent In a cosimulation context, the principle of linking the derivatives of the states to themselves (and potentially with time) is the same as in the monolithic case, yet the inputs and outputs also have to be considered.

Let $\nsys\in\mathbb{N}^*$ be the number of subsystems. Please note that the case $\nsys=1$ corresponds to a monolithic system. The cases that will be considered here are connected subsystems, that is to say $\nsys \geqslant 2$ subsystems that need to exchange data to one another: each input of each subsystem has to be fed by an output of another subsystem. A subsystem will be referenced by its subscript index $k\in\Insys$ so that subsystem-dependent functions or variables will have an index indicating the subsystem they are attached to.

Considering the inputs and the outputs, the co-simulation version of \eqref{eq:ODE} - \eqref{eq:ODE_decl} for subsystem $k\in\Insys$ is:
\vspace{-0.3cm}

\begin{equation}
\label{eq:ODE_eqonly_cosim}
\left\{
	\begin{array}{lll}
		\dsp{\frac{d}{dt}}x_k & = & f_k(t, x_k, u_k)\\
		y_k             & = & g_k(t, x_k, u_k)\\
	\end{array}
\right.
\end{equation}

\vspace{-0.3cm}
\noindent where
\vspace{-0.3cm}

\begin{equation}
\label{eq:ODE_eqonly_cosim_decl}
\begin{array}{l}
x_k \in \rnstk{k},\
u_k \in \rnink{k},\
y_k \in \rnoutk{k}\\
f_k:[\tinit, \tend] \times \rnstk{k} \times \rnink{k} \rightarrow \rnstk{k}\\
g_k:[\tinit, \tend] \times \rnstk{k} \times \rnink{k} \rightarrow \rnoutk{k}\\
\end{array}
\end{equation}

Equations \eqref{eq:ODE_eqonly_cosim} - \eqref{eq:ODE_eqonly_cosim_decl} are the equations representing a given subsystem. Ther are the minimal data required to entirely characterize any subsystem, yet they do not define the whole co-simulation configuration (connections are missing and can be represented by extra data, see \textit{link function} \cite{Eguillon2019}).

To stick to precise concepts, we should write $x$ as a function $x:[\tinit, \tend]\rightarrow\mathbb{R}^{n_{st}}$ (respectively for $y$, and $u$) as it is a time-dependent vector. That being said, we will keep on using $x$, $y$ and $u$ notations, as if they were simple vectors.

This representation enables us to evaluate derivatives at a precise point that has been reached. It is compliant to methods that need only the subsystems' equations, such as \textit{Decoupled Implicit Euler Method} and \textit{Decoupled Backward Differentiation Formulas} in \cite{Skelboe1992}.

\begin{figure}[!ht]
  \vspace{-0.2cm}
  \centering
   {\epsfig{file = 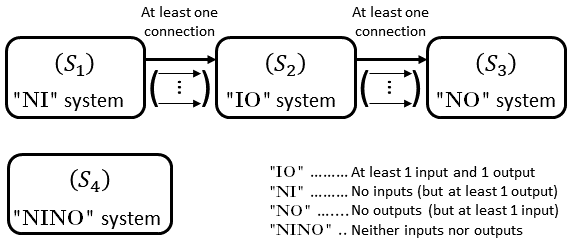, width = 7.0cm}}
  \caption{The four different topologies of subsystems, depending on the number of their coupling variables}
  \label{fig:NI_NO_NINO}
  \vspace{-0.1cm}
\end{figure}

When $\nink{k}$ and/or $\noutk{k}$ are $0$ for a given subsystem $k\in\Insys$, we will tag the \textit{topology} of the subsystem $(S_k)$ with a special name: NI, NO, NINO, or IO depending on the case (see figure \ref{fig:NI_NO_NINO}). This will be usefull to treat different behaviors when scheduling.

As a co-simulation implies communications at discrete times, let's introduce discrete notations. Let $n\in\mathbb{N}$ be the time index. Subsystem $(S_l)$ (with $l\in\Insys$) communicates at times $t_l^{[0]}, t_l^{[1]}, t_l^{[2]}, ...$. At these times, the $\noutk{l}$ outputs of $(S_l)$ are known: the value of the $j$\up{th} output (with $j\in\Inoutk{l}$) of subsystem $(S_l)$ at time $t_l^{[n]}$ will be written $y_{l, j}^{[n]}$.

Let $(S_k)$ be a subsystem with $\nink{k}>0$. The $i$\up{th} input is given by $u_{k, i}^{[n]}$ on macro-step $[t_k^{[n]}, t_k^{[n+1]}[$.

Please note that $\forall n\in\mathbb{N}$:

\noindent $\forall l\in\Insys, \forall j\in\Inoutk{l}, y_{l, j}^{[n]}\in\mathbb{R}$ whereas 

\noindent $\forall k\in\Insys, \forall i\in\Inink{k},$ $u_{k, i}^{[n]} : [t_k^{[n]}, t_k^{[n+1]}[ \rightarrow \mathbb{R}$.

\subsection{Polynomial calibration}
\label{subsection:polynomial_calibration}

\noindent For $q\in\mathbb{N}^*$, let $\mathbb{A}(q)\subsetneq \mathbb{R}^q$ be the set of elements that do not have any common values at different coordinates. In other words:
\vspace{-0.5cm}

\begin{equation}
\label{eq:distinguished_set}
\begin{array}{l}
	\mathbb{A}(q) = \big\{(t^{[r]})_{r\in[\![1, q]\!]}\in\mathbb{R}^q\ \big| \\
	\hspace{0.5cm} \forall (r_1, r_2)\in[\![1, q]\!]^2,\ r_1\neq r_2\Rightarrow t^{[r_1]}\neq t^{[r_2]}\big\}
\end{array}
\end{equation}
\vspace{-0.2cm}

For a given set of $q$ points $(t^{[r]}, z^{[r]})_{r\in[\![1, q]\!]}$ whose abscissas satisfy $(t^{[r]})_{r\in[\![1, q]\!]} \in \mathbb{A}(q)$, we define the two following polynomials:
\vspace{-0.3cm}

\begin{equation}
\label{eq:extrapolation_function}
\Omega_{q-1}^{Ex}:\left\{
\begin{array}{lcl}
	\mathbb{A}(q) \times \mathbb{R}^{q} \times \mathbb{R} & \rightarrow & \mathbb{R} \\
	\left(
		\begin{array}{c}
			(t^{[r]})_{r\in[\![1, q]\!]} \\
			(z^{[r]})_{r\in[\![1, q]\!]} \\
			t
		\end{array}
	\right)&
	\mapsto & \Omega_{q-1}^{Ex}\left(t\right) \\
\end{array}
\right.
\end{equation}
\vspace{-0.7cm}

\begin{equation}
\label{eq:constrainedleastsquares_function}
\Omega_{q-2}^{CLS}:\left\{
\begin{array}{lcl}
	\mathbb{A}(q) \times \mathbb{R}^{q} \times \mathbb{R} & \rightarrow & \mathbb{R} \\
	\left(
		\begin{array}{c}
			(t^{[r]})_{r\in[\![1, q]\!]} \\
			(z^{[r]})_{r\in[\![1, q]\!]} \\
			t
		\end{array}
	\right)&
	\mapsto & \Omega_{q-2}^{CLS}\left(t\right) \\
\end{array}
\right.
\end{equation}
\vspace{-0.3cm}

\noindent respectively called the \textit{extrapolation} and the \textit{constrained least squares} polynomials\footnote{For the sake of readability, we will sometimes write only the last variable of $\Omega_{q-1}^{Ex}$ and $\Omega_{q-2}^{CLS}$.}.

These polynomials are defined in order to have $\Omega_{q-1}^{Ex}\in\mathbb{R}_{q-1}[t]$ and $\Omega_{q-2}^{CLS}\in\mathbb{R}_{q-2}[t]$, where $\forall p\in\mathbb{N},\ \mathbb{R}_p[t]$ is the set of polynomials of the variable $t$ with coefficients in $\mathbb{R}$ and with a degree lower or equal to $p$. We have:
\vspace{-0.4cm}

\begin{equation}
\label{eq:extrapolation_function_def}
\forall r \in [\![1, q]\!], \Omega_{q-1}^{Ex}(t^{[r]}) = z^{[r]}
\end{equation}

\vspace{-0.3cm}
\noindent and
\vspace{-0.6cm}
	
\begin{equation}
\label{eq:constrainedleastsquares_function_def}
\begin{array}{c}
	\Omega_{q-2}^{CLS} : t \mapsto \sum_{i=0}^{q-2}a_i t^i\ \text{where} \\
	(a_i)_{\substack{\\ i\in[\![0, q-2]\!]}}
	=
	\!\!\!\!\!\!\!
	\underset
	{
		\substack
		{
			(\bar{a}_i)_{i\in[\![0, q-2]\!]} \\
			\Omega_{q-2}^{CLS} (t^{[1]}) = z^{[1]}
		}
	}
	{\arg \min \vspace{-0.1cm}}
	\!\!\!\!\!\!
	\Big\{
			\underset{r=1}{\overset{q}{\sum}}
			\!\!
			\Big(
				\!\!
				z^{[r]}
				\!
				-
				\!\!
				\underset{i=0}{\overset{q-2}{\sum}}
					\!
					\bar{a}_i (t^{[r]})^i
				\!
			\Big)^2
	\Big\}
\end{array}
\end{equation}
\vspace{-0.4cm}

In practice, the $z$ variables will be either inputs $u$ or outputs $y$, and index will correpond to a time index where $z^{[1]}$ is the latest one, so that the constraint on the least squares \eqref{eq:constrainedleastsquares_function_def} is the equality on the most recent point.

When we will consider one of these polynomials of degree $p$ (either extrapolation on $p+1$ points or constrained least squares on $p+2$ points) without specifying which one is used, we will use the generic notation $\Omega_p$.

The coefficients of $\Omega_p^{Ex}$  can be computed by several methods (Lagrange polynomials, Newton's formula, barycentric approach \cite{Berrut2004}), and the coefficients of $\Omega_p^{CLS}$ can be obtained with a cosntrained linear model (formula $(1.4.11)$ page $22$ of \cite{Amemiya1985}).

Finally, we define the Hermite interpolation polynomial $\mathcal{H}$ in the specific case with two points and first order derivatives:
\vspace{-0.5cm}

\begin{equation}
\mathcal{H}:\left\{
	\begin{array}{lcl}
		\mathbb{R}^2 \times \mathbb{R}^2 \times \mathbb{R}^2 \times \mathbb{R} &
		\rightarrow &
		\mathbb{R}
		\\
		(
			\left(
				\substack{t^{[1]}\\t^{[2]}}
			\right),
			\left(
				\substack{z^{[1]}\\z^{[2]}}
			\right),
			\left(
				\substack{\dot{z}^{[1]}\\ \dot{z}^{[2]}}
			\right),
			t
		) &
		\mapsto &
		\mathcal{H}(t)
	\end{array}
\right.
\end{equation}

\vspace{-0.2cm}
\noindent where
\vspace{-0.6cm}

\begin{equation}
\label{eq:hermite_function_def}
\begin{array}{l}
	\mathcal{H} \in \mathbb{R}_3[t]\ \text{and}\ \forall r \in [\![1, 2]\!], \\
	\mathcal{H}(t^{[r]}) = z^{[r]}
	\ \text{and}\ 
	\ndspfrac{d\mathcal{H}}{dt}(t^{[r]}) = \dot{z}^{[r]} \\
\end{array}
\end{equation}
\vspace{-0.2cm}

The coefficients of $\mathcal{H}$ can be computed with the square of the Lagrange polynomial basis or by using divided differences \cite{Hildebrand1956}.

Hermite interpolation will be used for smoothness enhancement in \ref{subsubsection:smoothness_enhancement}.

\section{\uppercase{\fornits algorithm}}
\label{sec:fornits_algorithm}

\noindent We introduce here the \fornits method, standing for Flexible Order Representation of New Inputs, including flexible Time-stepper (with variable step size, when applicable) and flexible Scheduler (asynchronous-capable, when applicable).

The method stems from the desire to keep the dynamical behavior of the coupling variables, what zero order hold (ZOH) does not do. At a given communication time, the outputs of the past data exchange will be reused in order to fit polynomial estimations for the future (the upcoming step). This is done in several co-simulation methods \cite{Kubler2000} \cite{Busch2016}, yet usually the polynomial order is decided in advance. \fornits has a flexible order and will decide of an order for every coupling variable at each communication time. Moreover, we will focus on the way to use these estimations in an asynchronous context (as it is not always possible to pass information at the same simulation time from one system to another). The error, depending on the polynomial order, will also be used to decide the evolution of the macro-step size (see figure \ref{fig:asynchronous_method}). Regarding subsystems with limited capabilities for holding time-dependent inputs, strategies are proposed to take advantage of the variable order thanks to an adaptation of the data to the capabilities of such subsystems. The latter strategies fit the specification to handle any modular model regardless the missing capabilities in every subsystem.

Finally, the time-stepping strategy (including the time-stepper itself and the scheduler) will be presented. The time-stepper is based on the error made by the estimation described in \ref{subsection:flexible_polynomial_inputs}. The normalization of this error has a strong impact on the time-stepping criterion, so several normalization methods will be described and a new one will be introduced: the \textit{damped amplitude} normalization method. The scheduler is also a part of the \fornits method, it occurs once the time-stepper produced an estimation of the upcoming step sizes, yet it shall not be fully detailed for the sake of space.

The smoothness enhancement \cite{Busch2019} \cite{Rauh2006} will also be presented as it is compliant with the \fornits method. Nevertheless, we adapted it to the context of a flexible order and variable step size method. The motivation is similar to the one in \cite{Eguillon2019}. In the case where subsystems do not have sufficient capabilities (up to $3$\up{rd} order polynomial inputs), the \fornits method can still run without smoothness enhancement.

\subsection{Flexible polynomial inputs}
\label{subsection:flexible_polynomial_inputs}

\noindent We arbitrary set the maximum degree for polynomial inputs to:
\vspace{-0.5cm}

\begin{equation}
\label{eq:max_degree}
M = 2
\end{equation}
\vspace{-0.5cm}

Let $(m_k)_{k\in\Insys}$ be the maximum degrees for polynomial inputs supported for each subsystem. As we want to support every kind of subsystem, we cannot assume anything on $m_k$ (we only know that $\forall k\in\Insys, m_k\geqslant 0$).

Let's define the effective maximum degree for each subsystem, by adding the constraint \eqref{eq:max_degree}:
\begin{equation}
\label{eq:max_degrees_k}
\forall k\in\Insys,\ M_k := \min(M, m_k)
\end{equation}

The time-dependent inputs that will be generated will always satisfy \eqref{eq:inputs_degrees_ranges}, and for each degree in these ranges, the maximum supported degree for each subsystem will never be exceeded thanks to \eqref{eq:max_degrees_k}.
\begin{equation}
\label{eq:inputs_degrees_ranges}
\forall k\in\Insys,\ \forall n\in\mathbb{N},\ u_k^{[n]} \in \left(\mathbb{R}_{M_k}[t]\right)^{\nink{k}}
\end{equation}

The determination of the order to use for polynomial inputs is quite straightforward. However, as this order may be different for each variable and as it is determined in the subsystem holding the corresponding variable as an output, some complications due to the asynchronousness may appear on subsystems having an input connected to this variable.

In order to clarify the process and to deal with properly defined mathematical concepts, we will split this explanation into three parts:

To begin with, we will define the \textit{order function} for every output variable and see how the values of this function are found. Then, we will define the estimated output variables based on order function. Finally, we will see how these estimated output variables are used from the connected input's perspective.

\subsubsection{Order function}
\label{subsubsection:order_function}

\noindent Let's consider a subsystem $(S_l)$ with $l\in\Insys$. Let's consider we have already done $n$ macro-steps with $n\in\mathbb{N}^*$, \ie at least one.

\begin{figure}[!ht]
  \vspace{-0.2cm}
  \centering
   {\epsfig{file = 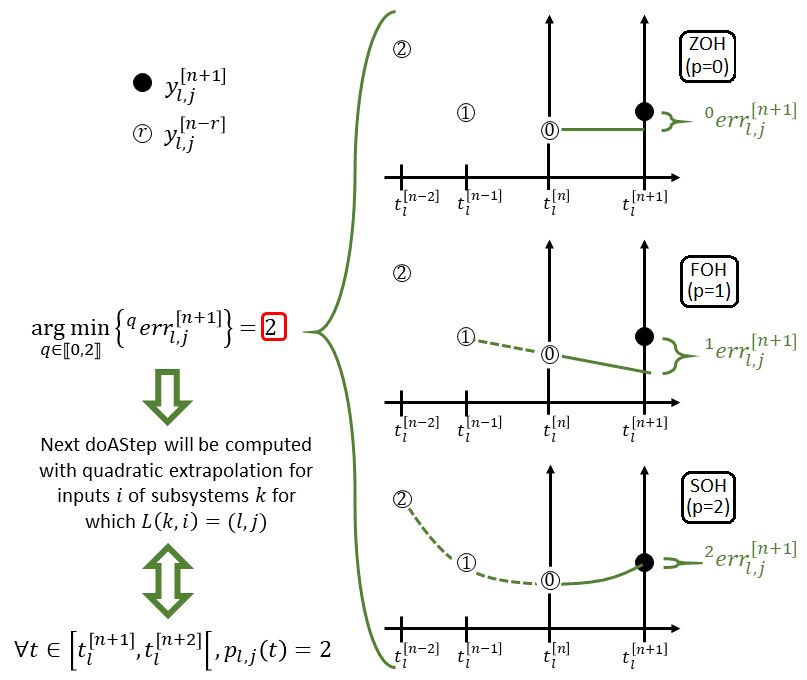, width = 7.5cm}}
  \caption{Order determination}
  \label{fig:order_determination}
  \vspace{-0.1cm}
\end{figure}

For $j\in\Inoutk{l}$ and $q\in[\![0, \min\left(M, n-1\right)]\!]$, we define $^{(q)}err_{l, j}^{[n+1]}$ the following way:
\begin{equation}
\label{eq:q_err_l_j_np1}
^{(q)}err_{l, j}^{[n+1]} = \left| y_{l, j}^{[n+1]} - \Omega_q^{Ex}\left(t_l^{[n+1]}\right) \right|
\end{equation}
\noindent where $\Omega_q^{Ex}$ is calibrated on $(t_l^{[n-r]})_{r\in[\![0, q]\!]}$ and $(y_{l, j}^{[n-r]})_{r\in[\![0, q]\!]}$ as explained in \eqref{eq:extrapolation_function} and \eqref{eq:extrapolation_function_def}. In other words the value obtained at $t_l^{[n+1]}$ is compared to the extrapolation estimation based on the $(q+1)$ exchanged values before $t_l^{[n+1]}$ (excluded), that is to say at $t_l^{[n-q]}, t_l^{[n-q+1]}, t_l^{[n-q+2]}, ..., t_l^{[n]}$ as shown on figure \ref{fig:order_determination}.

We now define the order function $p_l$ for subsystem $(S_l)$, which is a vectorial function for which each coordinate corresponds to one output of the subsystem $(S_l)$. It is a step function defining the "best order" of extrapolation for each macro-step, based on the previous macro-step. This includes a delay that makes it possible to integrate the subsystems simultaneously.

\begin{equation}
\label{eq:order_function_l}
\renewcommand{\arraystretch}{1.5}
p_l:\left\{
\begin{array}{lcl}
	[\tinit, \tend[ & \rightarrow & [\![0, M]\!]^{\noutk{l}} \\
	t & \mapsto	& \left(p_{l, j}(t)\right)_{j\in\Inoutk{l}} \\
\end{array}
\right.
\renewcommand{\arraystretch}{1}
\end{equation}

For $j\in\Inoutk{l}$, the $p_{l, j}$ functions are defined this way:

\begin{equation}
\label{eq:order_function_lj}
\renewcommand{\arraystretch}{1.5}
\begin{array}{l}
	p_{l, j}:\left\{
	\begin{array}{lcl}
		[\tinit, \tend[ & \rightarrow & [\![0, M]\!] \\
		t & \mapsto	& p_{l, j}(t) \ \text{with:}\\
	\end{array}
	\right. \\
	p_{l, j}(t) = 
		\dsp{\sum_{n=1}^{n_{\max}}}
		\left(
			\mathds{1}_{[t_l^{[n]}, t_l^{[n+1]}[}(t)\!\!\!\!\!\!\!
			\underset{q\in[\![0, \min(M, n-1)]\!]}{\arg\min}\!\!
			\left\{
				^q err_{l, j}^{[n]}
			\right\}
		\right)
\end{array}
\renewcommand{\arraystretch}{1}
\end{equation}

\noindent where $n_{\max}$ satisfies $t_l^{[n_{\max}]}=\tend$.

An illustration of the "$\arg\min$" choice in \eqref{eq:order_function_lj} is presented on figure \ref{fig:order_determination}, and the order function itself can be visualized on its plot in figure \ref{fig:order_function_p}.

\subsubsection{Estimated outputs}
\label{subsubsection:estimated_outputs}

\noindent We still consider the subsystem $(S_l)$ and its $\noutk{l}$ outputs. Let's admit the step $[t_l^{[n]}, t_l^{[n+1]}[$ is computed. Therefore, for all $j\in\Inoutk{l}$, the value of $p_{l, j}(t)$ for $t$ in $[t_l^{[n+1]}, t_l^{[n+2]}[$ is known.

We can now determine the \textit{estimated outputs} $\hat{y}_l^{[n+1]} \in\left(\mathbb{R}_{q}[t]\right)^{\noutk{l}},\ t\in [t_l^{[n+1]}, t_l^{[n+2]}[$. We have two choices for the way this estimation is made: \textbf{Extrapolation} mode and \textbf{Constrained Least Squares (CLS)} mode.

For the sake of genericity, we will define the abstraction of this choice $\Omega_q$ introduced in \ref{subsection:polynomial_calibration}.

Estimated outputs are defined the following way on the step $[t_l^{[n+1]}, t_l^{[n+2]}[$:

\begin{equation}
\label{eq:estimated_outputs_as_restricted_Omega_q}
\hat{y}_l^{[n+1]}:\left\{
\begin{array}{lcl}
	[t_l^{[n+1]}, t_l^{[n+2]}[ & \rightarrow & \rnoutk{l} \\
	t & \mapsto &
		\left(
			\hat{y}_{l, j}^{[n+1]}(t)
		\right)_{j\in\Inoutk{l}} \\
\end{array}
\right.
\end{equation}

\vspace{-0.2cm}
\noindent where
\vspace{-0.5cm}

\begin{equation}
\label{eq:estimated_output_j_as_restricted_Omega_q}
\hat{y}_{l, j}^{[n+1]}:\left\{
\begin{array}{lcl}
	[t_l^{[n+1]}, t_l^{[n+2]}[ & \rightarrow & \mathbb{R} \\
	t & \mapsto &
		\restrict
		{
			\Omega_{q}
		}
		{
			\big[
				t_l^{[n+1]}, t_l^{[n+2]}
			\big[
		}
		\left(
			t
		\right)\\
\end{array}
\right.
\end{equation}

where $\Omega_q$ is calibrated on $(t_l^{[n+1-r]}, y_{l, j}^{[n+1-r]})_{r\in[\![0, p]\!]}$ in extrapolation mode (see \eqref{eq:extrapolation_function_def}), on $(t_l^{[n+1-r]}, y_{l, j}^{[n+1-r]})_{r\in[\![0, p+1]\!]}$ in CLS mode (see \eqref{eq:constrainedleastsquares_function_def}), and where $q:=p_{l, j}(t_l^{[n+1]})$ so that $\deg(\hat{y}_{l, j}^{[n+1]}) = p_{l, j}(t_l^{[n+1]})$ (order determined previously is respected).

Let $\breve{y}_{l, j}^{[n+1]}$ be the extension of the polynomial $\hat{y}_{l, j}^{[n+1]}$ on the whole $\mathbb{R}$ domain.

\begin{figure}[!ht]
  \vspace{-0.2cm}
  \centering
   {\epsfig{file = 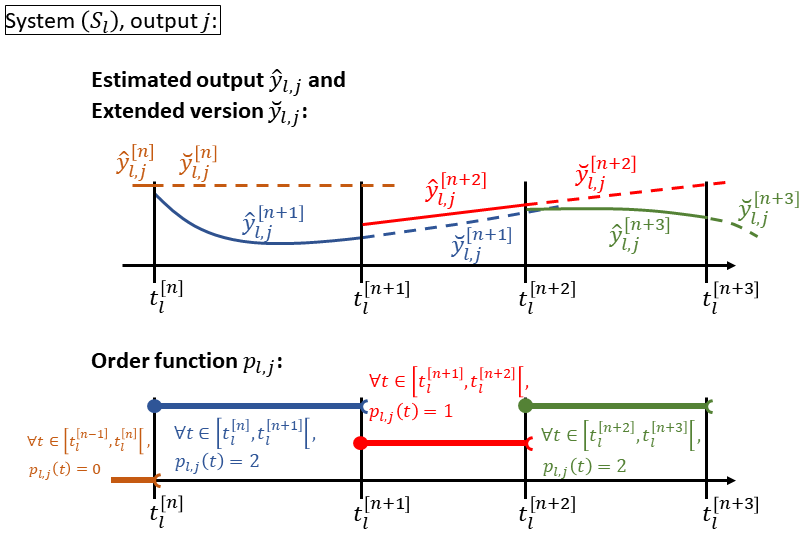, width = 7.5cm}}
  \caption{Order function and corresponding output estimation with corresponding orders}
  \label{fig:order_function_p}
  \vspace{-0.1cm}
\end{figure}

The degree of $\breve{y}_{l, j}^{[n+1]}$ is given by $p_{l, j}(t_l^{[n+1]})$ as shown on figure \ref{fig:order_function_p}.

As in \eqref{eq:estimated_output_j_as_restricted_Omega_q}, for the rest of \ref{subsubsection:estimated_outputs}, we will use the following notation (for the sake of readability):
\vspace{-0.2cm}

\begin{equation}
\label{eq:q_shortcut_np1}
q = p_{l, j}(t_l^{[n+1]})
\end{equation}

\begin{figure}[!ht]
  \vspace{-0.2cm}
  \centering
   {\epsfig{file = 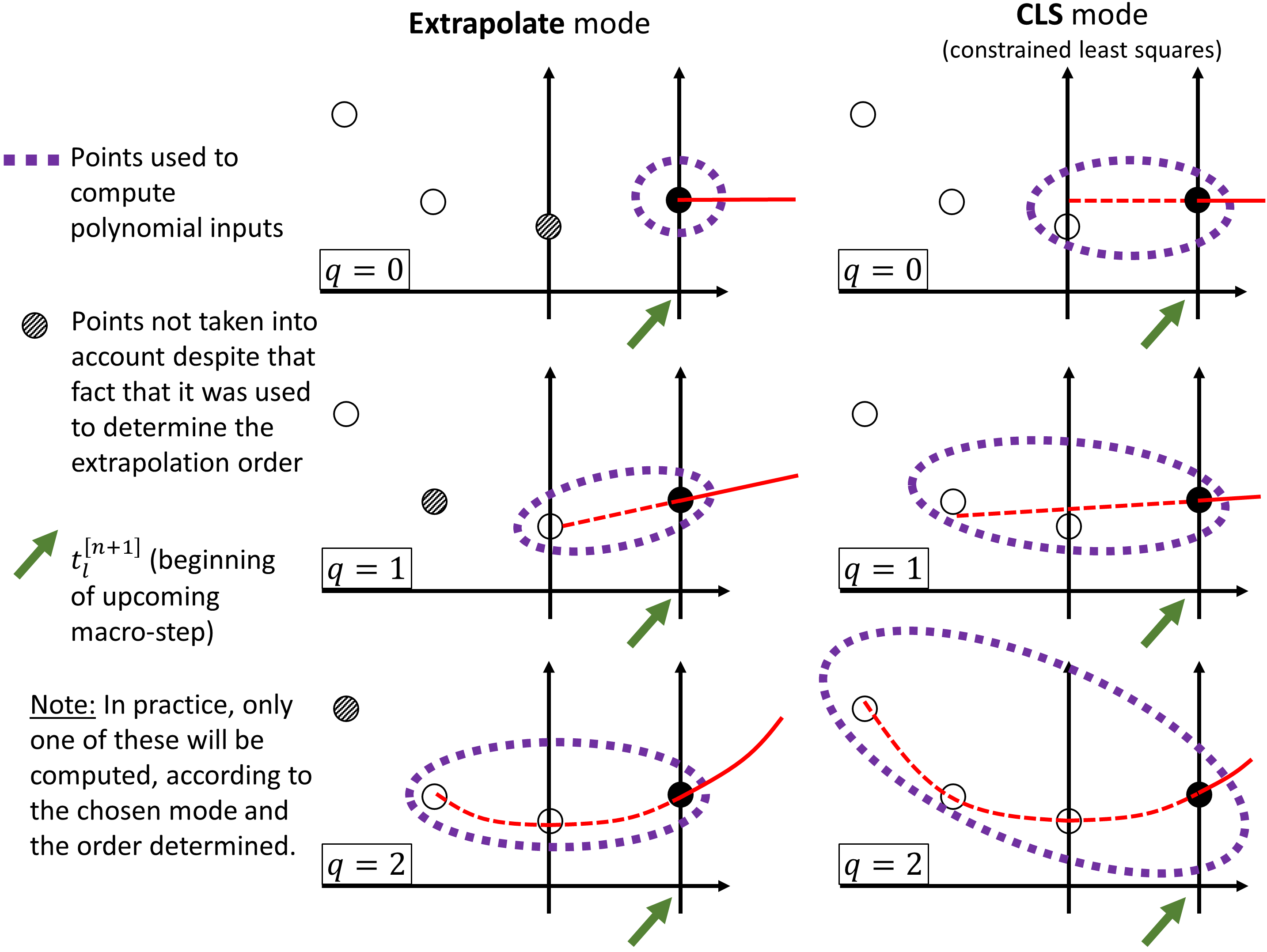, width = 7.5cm}}
  \caption{Extrapolation mode vs. CLS mode, on step $[t_l^{[n+1]}, t_l^{[n+2]}[$ once step $[t_l^{[n]}, t_l^{[n+1]}[$ is done (so that $p_{l, j}(t^{[n+1]})$ could be computed)}
  \label{fig:extrapolation_vs_constrainedleastsquares}
  \vspace{-0.1cm}
\end{figure}

The definition of $\Omega_q$ varies depending on the mode. In "Extrapolation" mode, the oldest point is forgotten (\ie not taken into account in the extrapolation): $(t_l^{[n-q]},\ y_{l, j}^{[n-q]})$. This point, taken into account to choose the order $q$ but not used in extrapolation calibration, is represented as striped in figure \ref{fig:extrapolation_vs_constrainedleastsquares}.

On the other hand, the idea of "CLS" mode is to take into account the point described above in the estimation of the output. An idea can be to forget the most recent point, that is to say: {\small $(t_l^{[n+1]}, y_{l, j}^{[n+1]})$}, but this would mean that the value given by the subsystem's integrator would be unused, so we will introduce a delay in the coupling process. Thus, the strategy is to take into account all the $(q\!+\!1)$ points that have been used in the determination of the chosen order, and the most recent point as well. We thus have $(q\!+\!2)$ points to adjust a polynomial of degree at most $q$: an extrapolation process cannot be made, but the "best fitting" polynomial can be found for the CLS criterion \eqref{eq:constrainedleastsquares_function_def}. Please note that removing the constraint $\Omega_{q-2}^{CLS}(t^{[1]})=z^{[1]}$ in \eqref{eq:constrainedleastsquares_function_def} corresponds to the relaxation technique on the past refered to as "method 1" in \cite{Li2020} in the particular case of $q=0$.

\subsubsection{Estimated inputs}
\label{subsubsection:estimated_inputs}

\noindent From the inputs perspective, the order function of the connected output is used. Let's consider a subsystem $(S_k)$ with $k\in\Insys$, which have $\nink{k}>0$ inputs. As the case $\nink{k}=0$ should not be excluded, let's admit that nothing is done from the inputs perspective in this case (because there is no input). From here and for the whole \ref{subsubsection:estimated_inputs} section, we will consider $\nink{k}\in\mathbb{N}^*$.

We will consider the input $i$ with $i\in\Inink{k}$ and, to properly consider the connected output, we stand $l\in\Insys$ and $j\in\Inoutk{l}$ so that input $i$ of subsystem $k$ is fed by output $j$ of subsystem $l$ (with the \textit{link function} notation of \cite{Eguillon2019}, we would write $(l, j) = L(k, i)$).

As asynchronousness should be supported, a special care should be made when the step $[t_k^{[n]}, t_k^{[n+1]}[$ does not fit into one single definition of $\hat{y}_{l, j}^{[m]}$ for an $m\in\mathbb{N}$. In this case, we will use:

\begin{equation}
\label{eq:estimated_inputs}
u_{k, i}^{[n]}:\left\{
\begin{array}{lcl}
	[t_k^{[n]}, t_k^{[n+1]}[ & \rightarrow & \mathbb{R} \\
	t & \mapsto & \breve{y}_{l, j}^{[m]} (t)
\end{array}
\right.
\end{equation}

where $\breve{y}$ is defined in \ref{subsubsection:estimated_outputs} and where

\begin{equation}
\label{eq:m}
m = \max\left\{
	m\in\mathbb{N}
	\ \big|\
	t_l^{[m]} \leqslant t_k^{[n]}
\right\}
\end{equation}

\v

\begin{figure}[!ht]
  \vspace{-0.2cm}
  \centering
   {\epsfig{file = 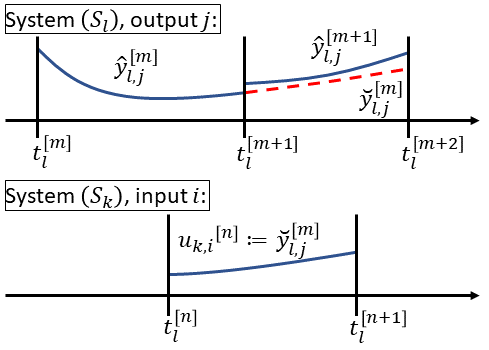, width = 6.5cm}}
  \caption{Estimated input using corresponding extended estimated output}
  \label{fig:extended_estimated_outputs}
  \vspace{-0.1cm}
\end{figure}

\eqref{eq:estimated_inputs} \eqref{eq:m} can be visualized figure \ref{fig:extended_estimated_outputs}.\\

Some subsystems cannot hold polynomial inputs, and some other can but not necessary at any order.

As \fornits method only requires to hold up to order $2$ (except when smoothness enhancement is triggered (see \ref{subsubsection:smoothness_enhancement}), which is a particular case), we will consider $9$ cases represented in the table below.

\begin{figure}[!ht]
  \vspace{-0.2cm}
  \centering
   {\epsfig{file = 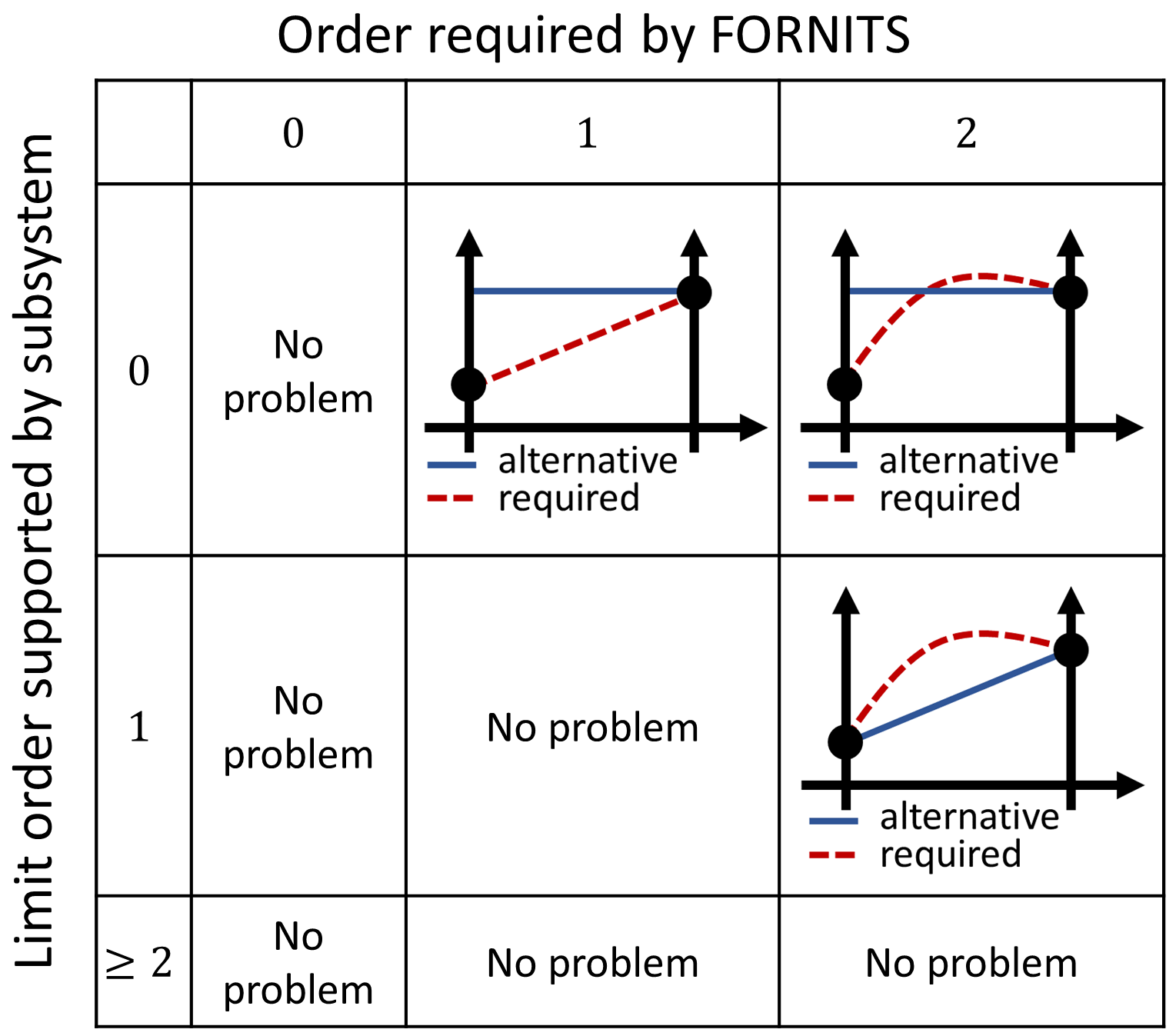, width = 6.5cm}}
  \caption{Alternatives to decrease polynomial input degree for subsystems with limited capabilities}
  \label{fig:table_of_degree_restriction}
  \vspace{-0.1cm}
\end{figure}

\subsubsection{Smoothness enhancement}
\label{subsubsection:smoothness_enhancement}

\noindent Smoothness enhancement can be triggered to enable $C^1$ inputs. This mode will only be applicable on subsystems supporting at least $3$\up{rd} order polynomial inputs.

Let's consider a system $(S_k)$ with $k\in\Insys$. Let's consider the step $[t_k^{[n]}, t_k^{[n+1]}[$ with $n\in\mathbb{N}^*$. We will also consider only the input $i$ of this system, with $i\in\Inink{k}$. As the smoothness enhancement process should be applied on every input separately, it will be detailed only on $u_{k, i}$ here.

The idea is to guarantee that $C^1$ smoothness is not broken at time $t^{[n]}$. In other word, we will remove the jump at the communication time. Moreover, as several consecutive steps are concerned, the $C^1$ smoothness won't be broken on the whole time interval (union of the steps).

\begin{figure}[!ht]
  \vspace{-0.2cm}
  \centering
   {\epsfig{file = 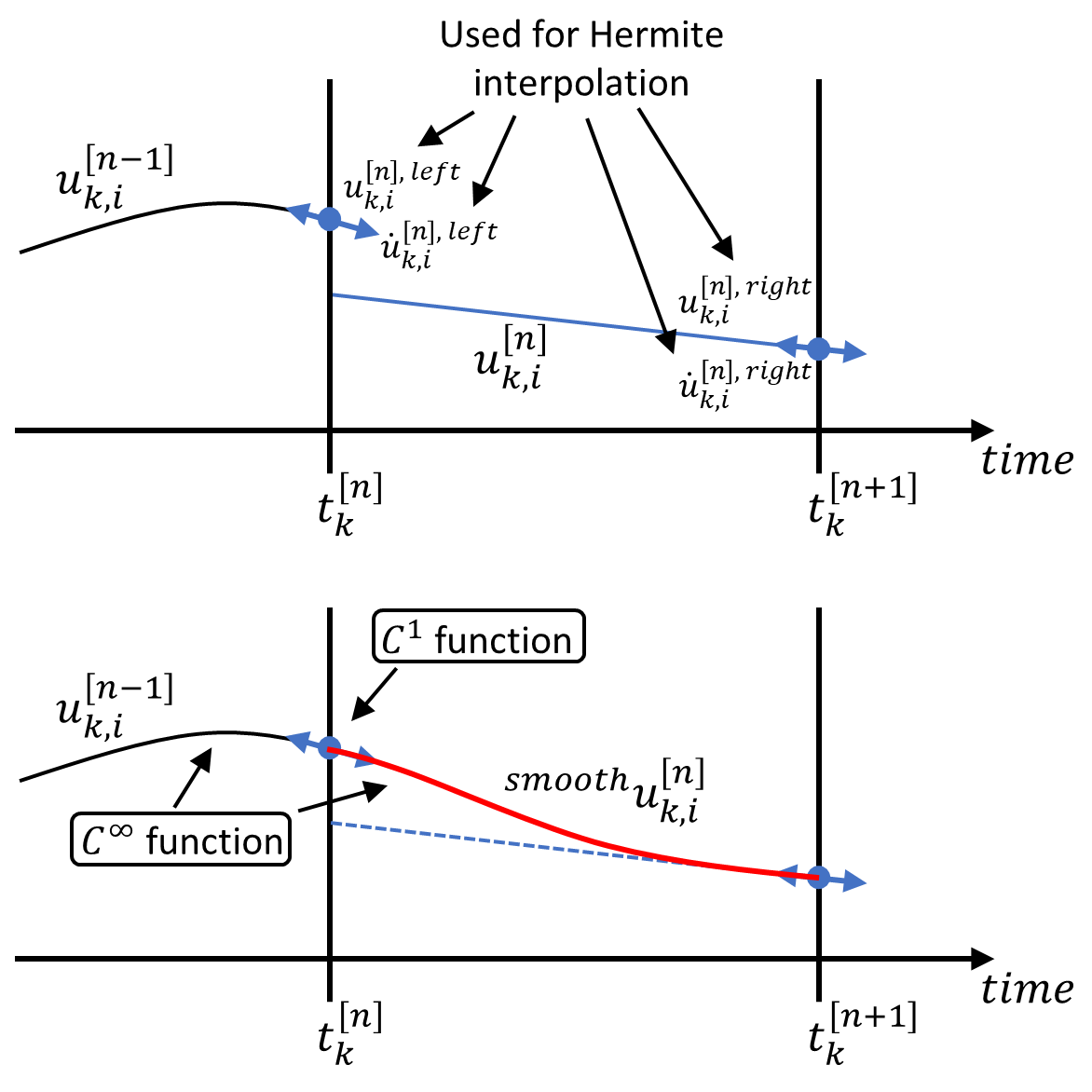, width = 7.5cm}}
  \caption{Redefinition of the time-dependent input on a co-simulation step in the case of smoothness enhancement}
  \label{fig:smoothness_enhancement}
  \vspace{-0.1cm}
\end{figure}

Regardless of the degree of the polynomial input, we will extend it to a third order polynomial using Hermite interpolation, as shown on figure \ref{fig:smoothness_enhancement}.

As $(S_k)$ reached the time $t_k^{[n]}$, we know the expression of the input $u_{k, i}^{[n-1]}$ that has been used for $[t^{[n-1]}, t^{[n]}[$, so we can compute the left constraints.
\vspace{-0.2cm}

\begin{equation}
\label{eq:C1_left_constraints}
\begin{array}{c}
	u^{[n], left} := \underset{\substack{t\rightarrow t^{[n]}\\t<t_{[n]}}}{\lim\ \ } u^{[n-1]}(t) \\
	\dot{u}^{[n], left} := \underset{\substack{t\rightarrow t^{[n]}\\t<t_{[n]}}}{\lim\ \ } \dot{u}^{[n-1]}(t)
\end{array}
\end{equation}
\vspace{-0.2cm}

Moreover, before performing the simulation on the step $[t^{[n]}, t^{[n+1]}[$, we compute the polynomial input (of degree $0$, $1$, or $2$, computed by extrapolation or by CLS) $u_{k, i}^{[n]}$. We thus compute the right constraints.
\vspace{-0.2cm}

\begin{equation}
\label{eq:C1_right_constraints}
\begin{array}{c}
	u^{[n], right} := u^{[n+1], left} := \underset{\substack{t\rightarrow t^{[n+1]}\\t<t_{[n+1]}}}{\lim\ \ } u^{[n]}(t) \\
	\dot{u}^{[n], right} := \dot{u}^{[n+1], left} := \underset{\substack{t\rightarrow t^{[n+1]}\\t<t_{[n+1]}}}{\lim\ \ } \dot{u}^{[n]}(t)
\end{array}
\end{equation}
\vspace{-0.2cm}

Finally, instead of using $u_{k, i}^{[n]}$ on $[t^{[n]}, t^{[n+1]}[$, we will use the "smooth" version of it:
\vspace{-0.2cm}

\begin{equation}
\label{eq:C1_inputs}
^{smooth}u_{k, i}^{[n]} : t \mapsto \mathcal{H}
\Big(
	\left(
		\substack
		{
		t^{[n]} \\
		t^{[n+1]}
		}
	\right),
	\left(
		\substack
		{
			u^{[n], left} \\
			u^{[n], right}
		}
	\right),
	\left(
		\substack
		{
			\dot{u}^{[n], left} \\
			\dot{u}^{[n], right}
		}
	\right),
	t
\Big)
\end{equation}

\noindent where $\mathcal{H}$ denotes the Hermite interpolation polynomial described in section \ref{subsection:polynomial_calibration}.

Using Hermite interpolation on values known by extrapolation is sometimes refered to as "extrapolated interpolation" \cite{Busch2019} \cite{Rauh2006}.

\subsection{Flexible time management}
\label{subsection:flexible_time_management}

\noindent The other major aspect of \fornits algorithm (the first one being the flexible polynomial inputs) is the time management, which includes a time stepper and a scheduler.

The time-stepper defines the next communication time after a macro-step is finished. The scheduler ensures the coherence of the rendez-vous times of all subsystems based on their connections, topologies and constraints.

\subsubsection{Time-stepper}
\label{subsubsection:time_stepper}

\noindent Let's consider a system $(S_l), l\in\Insys$ that is either IO or NI (see figure \ref{fig:NI_NO_NINO}). We have: $\noutk{l}>0$. Let's consider an output $y_{l, j}, j\in\Inoutk{l}$ of $(S_l)$.

The aim of the time-stepper is to determine $t_l^{[n+2]}$ once step $[t_l^{[n]}, t_l^{[n+1]}[$ have been computed.

For the sake of readability, we will use the notation $p:=p_{l, j}(t_l^{[n]})$ (different from $q$ in \eqref{eq:q_shortcut_np1} as here we focus on the lastly computed step, and not the upcoming one). Let's introduce the macro-step size and the dilatation coefficient, respectively \eqref{eq:macro_step_size} and \eqref{eq:dilatation_coefficient}.
\vspace{-0.2cm}

\begin{equation}
\label{eq:macro_step_size}
\dt_l^{[n]} = t_l^{[n+1]} - t_l^{[n]}
\end{equation}

\vspace{-0.5cm}

\begin{equation}
\label{eq:dilatation_coefficient}
\rho_l^{[n+1]} = \nicefrac{\dt_l^{[n+1]}}{\dt_l^{[n]}}
\end{equation}

They will be detailed later in this subsection, but for now we only need to know that the dilatation coefficients are bounded: $\exists (\rho_{\min}, \rho_{\max})\in\mathbb{R}_{+}^{*}, \forall l\in\Insys, \forall n\in\mathbb{N}, \rho_l^{[n]}\in[\rho_{\min}, \rho_{\max}]$.

In extrapolation mode, it exists $\zeta_t\in
\Big[
	t_l^{[n-p]},
	t_l^{[n+1]}
\Big]$ so that
\vspace{-0.2cm}

\begin{equation*}
\label{eq:extrapolation_error_applied_part_1}
\renewcommand{\arraystretch}{1.6}
\begin{array}{l}
	y_{l, j}^{[n+1]}-\breve{y}_{l, j}^{[n]}(t_l^{[n+1]}) \\
	=
		\underbrace{\dspfrac{1}{(p+1)!}\ddd{^{(p+1)}y_{l, j}}{t^{(p+1)}}(\zeta_t)}_{\text{independent of } \dt_l^{[n]}} \dspprod_{r=0}^{p}(t_l^{[n+1]}-t_l^{[n-r]}) \\
	=
		c_1\cdot \dspprod_{r=0}^{p}\Big(\sum_{s=0}^{r}\delta t_l^{[n-s]}\Big) \\
	=
		c_1\cdot \dspprod_{r=0}^{p}\Big(\dt_l^{[n]} \Big(1 + \sum_{s=1}^{r}\Big(\prod_{c=1}^{s}\ndspfrac{1}{\underbrace{\rho_l^{[n+1-c]}}_{\in[\rho_{\min}, \rho_{\max}]}}\Big)\Big)\Big) \\
\end{array}
\renewcommand{\arraystretch}{1.0}
\end{equation*}

\begin{equation}
\label{eq:extrapolation_error_applied}
\renewcommand{\arraystretch}{1.6}
\begin{array}{l}
	\leqslant
		c_1\cdot \dspprod_{r=0}^{p}\Big(\dt_l^{[n]} \Big(1 + \sum_{s=1}^{r}\Big(\prod_{c=1}^{s}\ndspfrac{1}{\rho_{\min}}\Big)\Big)\Big) \\
	\leqslant
		c_1\cdot \Big(\dt_l^{[n]}\Big)^{p+1} \underbrace{\dspprod_{r=0}^{p}\Big(\sum_{s=0}^{r}\Big(\ndspfrac{1}{\rho_{\min}}\Big)^{s}\Big)}_{\text{independent of } \dt_l^{[n]}} \\
	\leqslant
		c_1\cdot \Big(\dt_l^{[n]}\Big)^{p+1} \cdot c_2 \\
\end{array}
\renewcommand{\arraystretch}{1.0}
\end{equation}
\vspace{-0.2cm}

\noindent where $c_1$ and $c_2$ are random constants.

Final expression of error in \eqref{eq:extrapolation_error_applied} shows that the error is of order $(p+1)$ on the mean macro-step sizes. Analogously, we will consider that the error is of the same order error  in CLS mode: $(p+1)$ for a polynomial of degree $p$ (generated with a constrained least-square fitting on $(p+2)$ points).

The time-stepper uses this known error and error order to adapt the step size accordingly. We use a formula similar to the one in \cite{Schierz2012} (also mentioned in \cite{Gomes2018survey}) to define a dilatation coefficient candidate per output $j\in\Inoutk{l}$.
\vspace{-0.2cm}

\begin{equation}
\label{eq:dilatation_coefficient_per_output}
\rho_{l, j}^{[n+1]}:= \sqrt[p+1]{\frac{1}{error_{l, j}^{[n+1]}}}
\end{equation}
\vspace{-0.2cm}

In \eqref{eq:dilatation_coefficient_per_output}, the $error$ term is expected to be a relative error either relative to the values of the concerned variable and to a relative tolerance given which will determine the error threshold over which the step size is expected to decrease.

Given such a tolerance $\rtol$ and an absolute tolerance $\atol$, a first approach refered to as \textbf{Magnitude} relies on the order of magnitude of the variable at the moment of the communication. In that case, the error is defined as follow.
\vspace{-0.5cm}

\begin{equation}
\label{eq:error_magnitude}
^{Magn.}error_{l, j}^{[n+1]} := \frac{\left|y_{l, j}^{[n+1]} - \breve{y}_{l, j}^{[n]}(t_l^{[n+1]})\right|}{\atol + \rtol \cdot y_{l ,j}^{[n+1]}}
\end{equation}
\vspace{-0.2cm}

The problem in this approach is that when values of $y_{l ,j}$ are close to zero, the trend is to give a big error (and the step size will then be reduced). It is particularly problematic with variables with a great order of magnitude and periodically crossing the value zero, such as sinusoids.

An approach that might reduce this effect is to normalize the error according to the amplitude (observed since the beginning of the simulation) instead of the order of magnitude. This approach will be refered to as \textbf{Amplitude} and defines the error as follow.
\vspace{-0.2cm}

\begin{equation}
\label{eq:error_amplitude}
\begin{array}{l}
	^{Ampl.}error_{l, j}^{[n+1]} :=
	\\ \\
	\multicolumn{1}{r}
	{
		\dspfrac
		{
			\left|y_{l, j}^{[n+1]} - \breve{y}_{l, j}^{[n]}(t_l^{[n+1]})\right|
		}
		{
			\atol
			+
			\rtol
			\!
			\left(
				\!
				\underset{m\in[\![0, n+1]\!]}{\max}\!\!\left(y_{l, j}^{[m]}\right)
				-\!\!\!\!
				\underset{m\in[\![0, n+1]\!]}{\min}\!\!\left(y_{l, j}^{[m]}\right)
				\!\!
			\right)
		}
	}
\end{array}
\end{equation}
\vspace{-0.2cm}

The problem in this approach is that when the values of $y_{l, j}$ undergoes a great jump at one single moment (\textit{e.g.} at initialization), the trend will be to produce artificially small errors. Therefore, big step sizes will be produced by the time-stepper, and the accuracy may dwindle.

Solving this problem can be done by damping the amplitude, in order to progressively erase the effects of jumps while keeping the local amplitude of the variable $y_{l, j}$. The error produced with this principle will refer to the \textbf{Damped Amplitude} strategy, and will be defined as follow.
\vspace{-0.5cm}

\begin{equation}
\label{eq:error_damped_amplitude}
\begin{array}{l}
	\!\!\!	\!\!\!	\!\!\!
	^{\substack{Damped\\Ampl.}}error_{l, j}^{[n+1]} :=
	\\ \\
	\multicolumn{1}{r}
	{
		\hspace{10mm}
		\dspfrac
		{
			\left|y_{l, j}^{[n+1]} - \breve{y}_{l, j}^{[n]}(t_l^{[n+1]})\right|
		}
		{
			\atol
			+
			\rtol
			\!
			\left(
				\!
				^{\substack{damp\\max}}\!\!y_{l, j}^{[n+1]}
				-
				^{\substack{damp\\min}}\!\!y_{l, j}^{[n+1]}
				\!\!
			\right)
		}
	}
\end{array}
\end{equation}
\vspace{-0.2cm}

Expression \eqref{eq:error_damped_amplitude} refers to the damped minimal and damped maximal sequences which are recursively defined as follow.
\vspace{-0.5cm}

\begin{equation}
\label{eq:damped_min_max}
\renewcommand{\arraystretch}{1.5}
\left\{
\begin{array}{>{\hspace{-0.2cm}}l<{\hspace{-0.3cm}}c<{\hspace{-0.3cm}}l}
	\alpha^{[0]} & = & 0 \\
	^{\substack{damp\\max}}\!\!y_{l, j}^{[0]} & = & y_{l ,j}^{[0]} \\
	^{\substack{damp\\min}}\!\!y_{l, j}^{[0]} & = & y_{l ,j}^{[0]} \\
	\alpha^{[m]} & = & 
		^{\substack{damp\\max}}\!\!y_{l, j}^{[m]}
		-
		^{\substack{damp\\min}}\!\!y_{l, j}^{[m]} \\
	^{\substack{damp\\max}}\!\!y_{l, j}^{[m]}
	& = &
	\max\!\Big\{
		y_{l ,j}^{[m]},
		^{\substack{damp\\max}}\!\!y_{l, j}^{m-1}
		\!-
			\ndspfrac{\nu \cdot \delta t_l^{[m-1]}}{2}
			\cdot
			\alpha^{[m-1]}
	\Big\} \\
	^{\substack{damp\\min}}\!\!y_{l, j}^{[m]}
	& = &
	\min\!\Big\{
		y_{l ,j}^{[m]},
		^{\substack{damp\\min}}\!\!y_{l, j}^{m-1}
		\!+
			\ndspfrac{\nu \cdot \delta t_l^{[m-1]}}{2}
			\cdot
			\alpha^{[m-1]}
	\Big\} \\
\end{array}
\right.
\renewcommand{\arraystretch}{1.0}
\end{equation}

\noindent where $\forall n\in\mathbb{N},\ \delta t_l^{[n]} = t_l^{[n+1]} - t_l^{[n]}$ denotes the step size.

The damping coefficient $\nu\geqslant 0$ has to be defined. The greater it will be, the faster an event such as a jump will be "forgotten"; the smaller it will be, the closer $^{\substack{Damped\\Ampl.}}error_{l, j}^{[n+1]}$ will be to $^{Ampl.}error_{l, j}^{[n+1]}$.

We can now provide a proper definition of the subsystem's dilatation ratio mentioned in \eqref{eq:dilatation_coefficient} as the safest $\rho$ among the candidates \eqref{eq:dilatation_coefficient_per_output}.
\vspace{-0.2cm}

\begin{equation}
\label{eq:dilatation_coefficient_per_system}
\rho_l^{[n+1]} := \underset{j\in\Inoutk{l}}{\min}\left(\rho_{l, j}^{[n+1]}\right)
\end{equation}
\vspace{-0.2cm}

With this dilatation ratio, we obtain an estimation of the next communication time for subsystem $(S_l)$:
\vspace{-0.2cm}

\begin{equation}
\label{eq:dilatation_ratio_application}
\begin{array}{lcl}
	\delta t_l^{[n]} & = & t_l^{[n+1]} - t_l^{[n]} \\
	\widetilde{\delta t}_l^{[n+1]} & = & \rho_l^{[n+1]} \cdot \delta t_l^{[n]} \\
	\tilde{t}_l^{[n+2]} & = & t_l^{[n+1]} + \widetilde{\delta t}_l^{[n+1]}
\end{array}
\end{equation}
\vspace{-0.2cm}

This next estimation temporary solves the problematic introduced at the beginning of \ref{subsubsection:time_stepper}, the determination of $t_l^{[n+2]}$ once $[t_l^{[n]}, t_l^{[n+1]}[$ is computed. Nonetheless, this is only an estimation as the scheduler presented in \ref{subsubsection:scheduler} may modify it and determine the $t_l^{[n+2]}$ to use.

\textbf{\underline{Remark 1 (ratio bounds):}} In practice, a minimal and a maximal value for the dilatation ratio \eqref{eq:dilatation_coefficient_per_system} can be define in order to avoid unsafe extreme step size reduction/increase. In section \ref{sec:results_and_behavior_on_two_test_cases}, the dilatation coefficient is projected in the interval $[10\%, 105\%]$.

\textbf{\underline{Remark 2 (subsystems without outputs):}} The above procedure only works for subsystems having at least one output variable. For the other subsystems (NO and NINO), we will use the following rule.
\vspace{-0.3cm}

\begin{equation}
\label{eq:time_stepper_NO_NINO}
\begin{array}{l}
	\forall l\in\left\{l\in\Insys\ |\ \noutk{l} = 0\right\}, \\
	\hspace{2cm}
		\forall n \in\mathbb{N}\cap[2, +\infty[,\ \tilde{t}_l^{[n]} = \tend
\end{array}
\end{equation}
\vspace{-0.2cm}

\textbf{\underline{Remark 3 (co-simulation start):}} The rules presented above do not enable to set the two first times. We will obviously set the first one with $\forall l\in\Insys, \tilde{t}_l^{[0]} = \tinit$ and the next one by using an initial step $\delta t_l^{[0]}$ given as co-simulation parameter for every subsystem. This initial step size will have a limited influence on the whole co-simulation. We then set $\tilde{t}_l^{[1]} = \tinit + \delta t_l^{[0]}$ for each subsystem.

\subsubsection{Scheduler}
\label{subsubsection:scheduler}

\noindent The scheduler is intended to adjust the communication times of the subsystems according to the times they reached individually. It is based on topological data (an adjacency matrix representing the connections between the subsystems and the set of constraints regarding step sizes on each subsystem) and can handle asynchronousness in a way that makes it robust to co-simulation involving subsystems with imposed step sizes. Several stages happen sequentially once every subsystem has produced an estimated next communication time (see \ref{subsubsection:time_stepper}). Among them, some are dedicated to avoid to use a variable further than the time it is supposed to be used according to the time-stepper (trying to avoid phenomena like the one happening on figure \ref{fig:extended_estimated_outputs}), and others act as optimisations, as they avoid to communicate data when this has no effect (NI systems feeding a NO subsystem with a large imposed step size, for instance).

The stages of the scheduler will not be detailed here for the sake of space. It can be seen as a black box producing the effective next communication steps once the time-stepper produced the estimated ones.

In other words, once $\forall l\in\Insys,\ [t_l^{[n_l]}, t_l^{[n_l+1]}[$ is computed ($n_l$ might be different per subsystems as some subsystems may be idle when other keep on going, so they might do a different amount of steps), the time-stepper has produces $(\tilde{l}^{[n_l+1]})_{l\in\Insys}$, and the scheduler can be seen as the black box $Sch$ below.
\vspace{-0.2cm}

\begin{equation}
\label{eq:scheduler}
(t_l^{[n_l+2]})_{l\in\Insys} := Sch
\left(
\begin{array}{c}
	(t_l^{[n_l+1]})_{l\in\Insys} \\
	(\tilde{t}_l^{[n_l+2]})_{l\in\Insys}
\end{array}
\right)
\end{equation}
\vspace{-0.2cm}

\section{\uppercase{Results and behavior on two test cases}}
\label{sec:results_and_behavior_on_two_test_cases}

\noindent Two test cases will be presented in this section. Both of them have been made with Simcenter Amesim.

The first model will present the importance that keeping the dynamics on the coupling variables may have. This is actually done with the polynomial inputs in \fornits method. The second model is a variation of the 2-masses-springs-dampers \cite{Eguillon2019} \cite{Busch2016} and performances will be measured on several co-simulations.

The most widely used co-simulation method is the non-iterative Jacobi (parallel) fixed-step algorithm, with zero-order-hold inputs. Therefore, the comparisons will be made between the latter (refered to as "NI Jacobi") and \fornits algorithm.

\subsection{Car with controlled speed}
\label{subsection:car_with_controlled_speed}

\noindent The model consists in two subsystems, one corresponds to a $1000$kg car (simplified as a mass) moving on a $1$D axis (straight road) modelled with Newton's second law, which takes a force on entry and gives its position as output, and the other corresponds to a controller, producing a force and using the car position on input.

\begin{figure}[!ht]
  \vspace{-0.2cm}
  \centering
   {\epsfig{file = 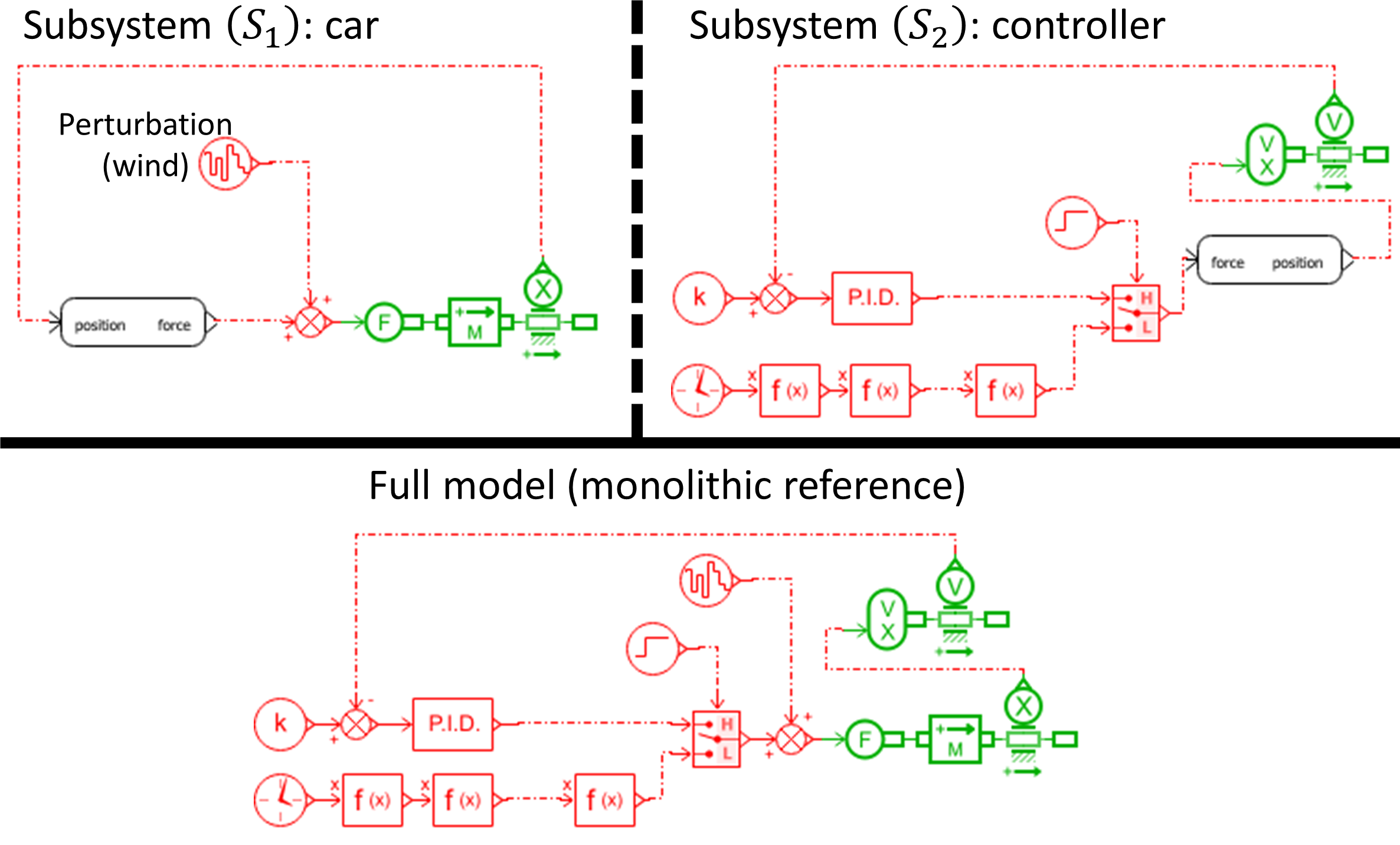, width = 7.5cm}}
  \caption{Test case 1: car with controlled speed - Subsystems and monolithic reference in Simcenter Amesim}
  \label{fig:Testcase_1}
  \vspace{-0.1cm}
\end{figure}

\begin{figure}[!ht]
  \vspace{-0.2cm}
  \centering
   {\epsfig{file = 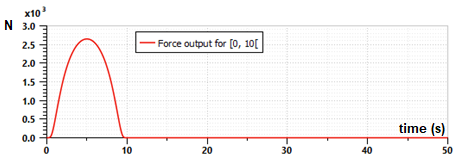, width = 6.0cm}}
  \caption{Preset output force from the controller on $[0, 10[$}
  \label{fig:Testcase_1_initial_force}
  \vspace{-0.1cm}
\end{figure}

As it can be seen on figure \ref{fig:Testcase_1}, the car subsystem adds a random perturbation to the input force (this may be seen as a 1D wind, for instance). Moreover, on the first $10$s of the (co-)simulation, the controller will output a force that is predetermined and that can be seen on figure \ref{fig:Testcase_1_initial_force}.

At $t=10$s, the output of the controller subsystem becomes the output of a "P" controller based on the vehicle speed and which is designed to make the vehicle reach (or maintain) a target speed of $16$m s\up{$-1$}. The velocity needed by the "P" controller is computed from the input position, using an explicit equation $\dot{x} = v$ and a constraint equation $x-u=0$ where $v$ denotes the computed speed, $x$ a hidden variable representing the position, and $u$ the position input in the controller subsystem.

\begin{figure}[!ht]
  \vspace{-0.2cm}
  \centering
   {\epsfig{file = 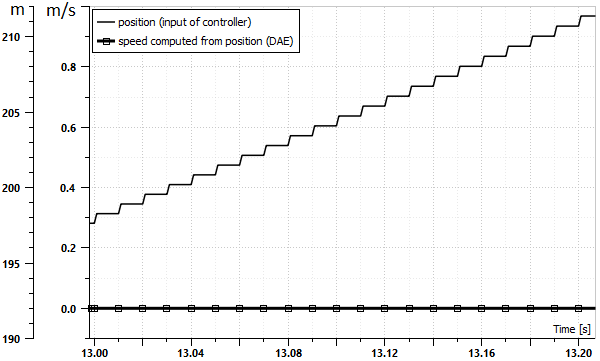, width = 6.0cm}}
  \caption{Vehicle speed computed in controller $(S_2)$ using its input (vehicle position) with a ZOH input}
  \label{fig:vehicle_speed_computed_in_controller_NI_Jacobi}
  \vspace{-0.1cm}
\end{figure}

Due to the zero order hold, the NI Jacobi method does not allow $(S_2)$ to properly retrieve the vehicle speed. Indeed, on every co-simulation step, the input position will be constant and the DAE will produce a null speed (see figure \ref{fig:vehicle_speed_computed_in_controller_NI_Jacobi}, as it can be seen on the zoom on $[13, 13.2]$.

\begin{figure}[!ht]
  \vspace{-0.2cm}
  \centering
   {\epsfig{file = 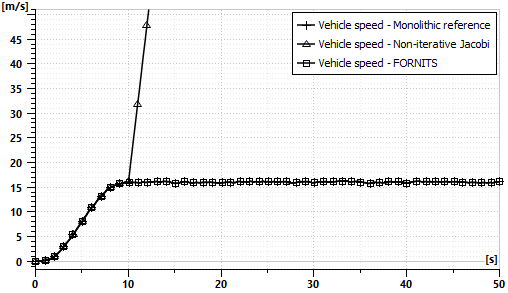, width = 6.0cm}}
  \caption{Vehicle speed in $(S_1)$ in the monolithic reference and in co-simulations with respectively the NI Jacobian (ZOH) method and the \fornits algorithm}
  \label{fig:vehicle_speed}
  \vspace{-0.1cm}
\end{figure}

The consequence of this is that the controller will put force to push the vehicle in order to increase the speed (attempting to make it reach $16$m s\up{$-1$}) without even realizing that the vehicle has a non-null velocity. Therefore, the velocity will keep on increasing indefinitely (curve with triangles on figure \ref{fig:vehicle_speed}). Inversely, \fornits represents the dynamics of the coupling variables (Force and vehicle position) when needed (this is the "flexible order" part presented in \ref{subsection:flexible_polynomial_inputs}). The velocity computed from such a position input is more relyable, and the controller $(S_2)$ will send proper force to make the vehicle maintain its target velocity.

The effective vehicle velocity inside of $(S_1)$ across the time is presented in figure \ref{fig:vehicle_speed} in the three cases of (co-)simulation.

\subsection{Two masses, springs and dampers}
\label{subsection:two_masses_springs_and_dampers}

\begin{figure}[!ht]
  \vspace{-0.2cm}
  \centering
   {\epsfig{file = 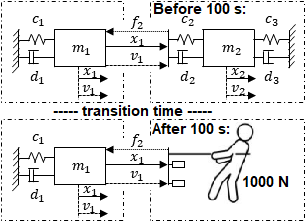, width = 7.2cm}}
   \vspace*{-0.2cm}
\begin{equation}
\label{eq:2_masses_params}
\renewcommand{\arraycolsep}{2.pt}
\begin{small}
\begin{array}{|lrr|lrr|lrr|}
	\hline
	c_1 & = &   1\ \text{kN/m} & \
		d_1 & = &   1\ \text{kN/(m/s)} & \
			x_1(\tinit) & = & -1\ \text{m} \\
	c_2 & = &   1\ \text{kN/m} & \
		d_2 & = &   0\ \text{kN/(m/s)} & \
			x_2(\tinit) & = & 0\ \text{m} \\
	c_3 & = &   1\ \text{kN/m} & \
		d_3 & = &   1\ \text{kN/(m/s)} &
			v_1(\tinit) & = & 0\ \text{m/s} \\
	m_1 & = & 1000\ \text{kg} & \
		m_2 & = & 1000\ \text{kg} &
			v_2(\tinit) & = & 0\ \text{m/s} \\
	\multicolumn{9}{|c|}{[\tinit, \tend] = [0\ \text{s},\  200\ \text{s}]} \\
	\hline
\end{array}
\end{small}\nonumber
\end{equation}
  \caption{Two masses model with coupling on force, displacement and velocity. Behavior of the right mass changes after $100$s}
  \label{fig:2_masses_cosim}
  \vspace{-0.1cm}
\end{figure}

The model consists of two masses coupled with force, displacement and velocity as shown on figure \ref{fig:2_masses_cosim}. It is a variant of the test cases presented in \cite{Eguillon2019} or in \cite{Busch2016} in the sense that, after a transition time (set to the middle of the simulation time interval), the behavior of one of the models changes. This generates a discontinuity on coupling variables and different behaviors before and after.

First of all, for the sake of performance and accuracy comparisons, the model have been co-simulated using the NI Jacobi method with five fixed step size values: $0.01$s, $0.05$s, $0.1$s, $0.2$s and $0.4$s. Then, \fornits method was used with $\dt_1^{[0]}=\dt_2^{[0]}=0.01$s (also used as "minimal step size value"), for different strategies regarding polynomial calibration (extrapolation / CLS, see \ref{subsection:polynomial_calibration}), smoothness enhancement (disabled / make interfaces $C^1$, see \ref{subsubsection:smoothness_enhancement}) and error normalization (w.r.t. amplitude / magnitude / damped amplitude with $\nu=5\%$, see \ref{subsubsection:time_stepper}).

\begin{figure}[!ht]
  \vspace{-0.2cm}
  \centering
   {\epsfig{file = 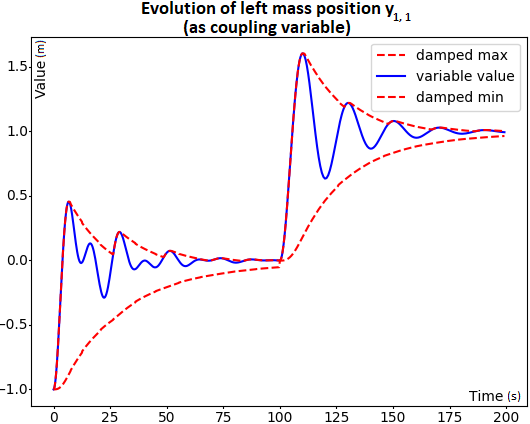, width = 5.0cm}}
  \caption{Evolution of $y_{1, 1}$ (coupling variable corresponding to the position of the left mass $x_1$) and its damped bounds sequences a defined in \eqref{eq:damped_min_max}}
  \label{fig:damped_amplitude}
  \vspace{-0.1cm}
\end{figure}


\begin{figure}[!ht]
  \vspace{-0.2cm}
  \centering
   {\epsfig{file = 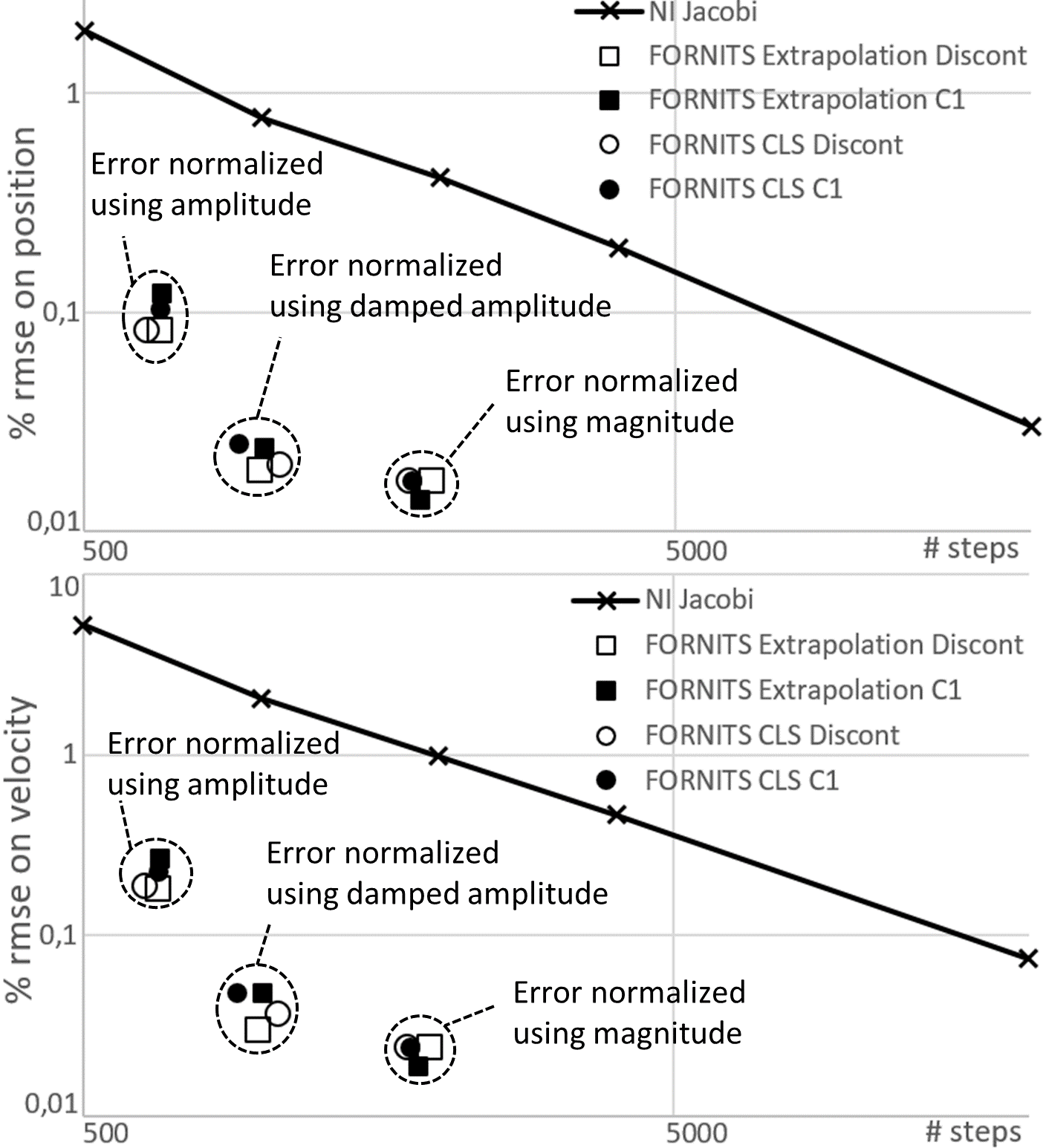, width = 7.0cm}}
  \caption{RMSE/\#steps trade-off (on $x_1$ and $v_1$) for various sets of parameters of \fornits algorithm on 2 masses model}
  \label{fig:2_masses_cosim_graphs}
  \vspace{-0.1cm}
\end{figure}

All co-simulations run instantly (measured elapsed time: $0.0$s), but a good performance indicator can be the number of steps proceeded. Indeed, the main cost of a co-simulation run is the frequent restart of subsystems' solvers after each discontinuity introduced by a communication time. On the other side, bigger co-simulation steps are expected to produce inaccurate results: the larger a co-simulation step is, the older an input value might be used. The precision criterion will be the RMSE (root mean square error) compared to the monolithic reference.

The error normalization criterion in the time-stepper has a strong impact on the results: as expected, $^{Magn.}error$ generates larger errors and smaller co-simulation steps. The total amount of steps is bigger than the other normalization methods, but the results are also more accurate. The trade-off is: results are either more accurate than with NI Jacobi for the same total number of steps, or obtained in less steps than NI Jacobi to achieve the same accuracy (see right-hand dashed circled sets of points on figure \ref{fig:2_masses_cosim_graphs}).

In order to reduce the effect of artificial step size reduction when coupling variables cross zero, the error normalization $^{Ampl.}error$ method may be used in the time-stepper. The steps are larger, and the total numbers of steps are the smallest ones on every case (compare the $3$\up{rd} subtable of table \ref{tab:results} to the $2$\up{nd} and the $4$\up{th} ones). However, as the model is damped ($d_1=d_3>0$), the initial highest and lowest values of the variables are not forgotten throughout the whole co-simulation (this can be observed on the blue curve of figure \ref{fig:damped_amplitude}). The trade-off is still better than NI Jacobi on every case.

In order to improve the accuracy of the results obtained using the $ ^{Ampl.}error$, the damped amplitude alternative for error normalization has been proposed in subsubsection \ref{subsubsection:time_stepper}. Indeed, damping the amplitude with a factor $\nu=5\%$ allows the coupling to progressively forget the initial high values, as the minimum and maximum values will "follow" the order of magnitude of the variables amplitude. The $^{\substack{damp\\max}}\!\!y_{1, 1}^{[n]}$ and $^{\substack{damp\\min}}\!\!y_{1, 1}^{[n]}$ sequences can be seen on figure \ref{fig:damped_amplitude}, $y_{1, 1}$ being the coupling variable associated to $x_1$.

Table \ref{tab:results} compiles the results of figure \ref{fig:2_masses_cosim_graphs} and shows precisely the rmse values. Among others, it shows that the polynomial inputs calibration method (either with an extrapolation $\Omega_q^{Ex}$ or a CLS fitting $\Omega_q^{CLS}$) has minor incidence on the final results on the two masses test case. Nonetheless, it can be observed on this table that the smoothness enhancement (presented in \ref{subsubsection:smoothness_enhancement}) slightly damages accuracy, except when the error is normalized using $ ^{Magn.}error$. The reason of this is the fact that the value of the inputs is not the one given as output by its connected subsystem at the beginning of the step. In other word, this phenomena (that can be observed at $t_k^{[n]}$ on figure \ref{fig:smoothness_enhancement}) generates a small extra error at the beginning of every co-simulation step, but the value of the input is quickly adapted to match the one it would have had without the smoothness enhancement. The smaller the steps are, the less this phenomena can be observed: that is the reason why the subtable regarding $ ^{Magn.}error$ does not show it. The benefit of smoothness enhancement is mainly the possibility to restart faster the solvers inside of the subsystems due to the $C^1$ smoothness of the input variables. The consequence of it would be a faster co-simulation run for the same number of co-simulation steps.

\vspace{1cm}

\begin{table}[!h]
\caption{Results on 2 masses test model - comparing number of co-simulation steps and rmse on $2$ state variables}\label{tab:results} \centering
\renewcommand{\arraystretch}{1.0}
\begin{tabular}{|c|c||c|c|c|}
\hline
	\multicolumn{5}{|c|}{NI Jacobi (ZOH)} \\
\hline
	\multicolumn{2}{|c||}{$\dt$} & \#steps & $\substack{\text{rmse}\\\text{on}\ x_1}$ & $\substack{\text{rmse}\\\text{on}\ v_1}$ \\
\hline
\hline
	\multicolumn{2}{|c||}{$0.01$} & $20 000$ & \footnotesize $0.030\%$ & \footnotesize $0.074\%$ \\
\hline
	\multicolumn{2}{|c||}{$0.05$} & $4 000$ & \footnotesize $0.197\%$ & \footnotesize $0.462\%$ \\
\hline
	\multicolumn{2}{|c||}{$0.1$} & $2 000$ & \footnotesize $0.412\%$ & \footnotesize $0.981\%$ \\
\hline
	\multicolumn{2}{|c||}{$0.2$} & $1 000$ & \footnotesize $0.772\%$ & \footnotesize $2.044\%$ \\
\hline
	\multicolumn{2}{|c||}{$0.4$} & $500$ & \footnotesize $1.928\%$ & \footnotesize $5.194\%$ \\
\hline
\multicolumn{5}{c}{ }\\
\hline
	\multicolumn{5}{|c|}{\fornits} \\
\hline
	$\!\Omega_q\!$ & {\footnotesize Smoothness} & \#steps & $\substack{\text{rmse}\\\text{on}\ x_1}$ & $\substack{\text{rmse}\\\text{on}\ v_1}$ \\
\hline
\hline
	\multicolumn{5}{|c|}{time-stepper using $^{Magn.}error$} \\
\hline
	\multirow{2}{*}{$\!\Omega^{Ex}_q\!$} & disabled & $1 935$ & \footnotesize $0.017\%$ & \footnotesize $0.024\%$ \\
\cline{2-5}
	& enhanced & $1 848$ & \footnotesize $0.014\%$ & \footnotesize $0.019\%$ \\
\hline
	\multirow{2}{*}{$\!\Omega^{CLS}_q\!$} & disabled & $1 768$ & \footnotesize $0.017\%$ & \footnotesize $0.024\%$ \\
\cline{2-5}
	& enhanced & $1 786$ & \footnotesize $0.017\%$ & \footnotesize $0.024\%$ \\
\hline
\hline
	\multicolumn{5}{|c|}{time-stepper using $^{Ampl.}error$} \\
\hline
	\multirow{2}{*}{$\!\Omega^{Ex}_q\!$} & disabled & $669$ & \footnotesize $0.083\%$ & \footnotesize $0.183\%$ \\
\cline{2-5}
	& enhanced & $674$ & \footnotesize $0.122\%$ & \footnotesize $0.269\%$ \\
\hline
	\multirow{2}{*}{$\!\Omega^{CLS}_q\!$} & disabled & $636$ & \footnotesize $0.083\%$ & \footnotesize $0.188\%$ \\
\cline{2-5}
	& enhanced & $672$ & \footnotesize $0.104\%$ & \footnotesize $0.227\%$ \\
\hline
\hline
	\multicolumn{5}{|c|}{time-stepper using damped amplitude error} \\
\hline
	\multirow{2}{*}{$\!\Omega^{Ex}_q\!$} & disabled & $989$ & \footnotesize $0.019\%$ & \footnotesize $0.030\%$ \\
\cline{2-5}
	& enhanced & $1007$ & \footnotesize $0.024\%$ & \footnotesize $0.048\%$ \\
\hline
	\multirow{2}{*}{$\!\Omega^{CLS}_q\!$} & disabled & $1071$ & \footnotesize $0.020\%$ & \footnotesize $0.037\%$ \\
\cline{2-5}
	& enhanced & $911$ & \footnotesize $0.025\%$ & \footnotesize $0.048\%$ \\
\hline
\end{tabular}
\renewcommand{\arraystretch}{1.0}
\end{table}

\section{\uppercase{Conclusions}}
\label{sec:conclusions}

\noindent The \fornits coupling algorithm used with a time-stepping criterion based on the local error estimation normalized with regard to magnitude gives a better error/\#steps trade-off that the non-iterative Jacobi method. The number of co-simulation steps can even be reduced by tuning the error normalization method, and the involved coefficient (refered to as $\nu$ damping coefficient). In other words, a safe approach could be to start with a co-simulation with \fornits using the error normalized with regard to the order of magnitude and, if the total number of steps is not satisfactory, the strategy can be to use the damped amplitude with a decent $\nu$ coefficient and to progressively decrease it ($\nu=0\%$ corresponds to the classical \textit{amplitude} approach).

On the second test case, with error normalized regarding amplitude damped with a coefficient $\nu=5\%$ (shapes in the middle dashed circles on figure \ref{fig:2_masses_cosim_graphs}), we can achieve, with approximately $20$ times less co-simulation steps (and as many avoided discontinuities), an accuracy similar to the one obtained with the non-iterative Jacobi method (or ever a bit more accurate results, as non-iterative Jacobi method does not reach an rmse lower than $0.03\%$ for $x_1$ even with $20\ 000$ steps). From another point of view, a similar amount of communication times produces an error $38$ times smaller with \fornits than with the non-iterative Jacobi method. Moreover, a co-simulation step size does not need to be chosen in advance as it is automatically adapted with \fornits on the subsystems allowing it.

The robustness of the method comes from the motivation at the starting point: \fornits can be applied regardless of the subsystems and their capabilities. It is therefore possible to integrate the method on an industrial product that must handle subsystems coming from a wide range of platforms.

The CLS approach, designed in order to avoid forgetting data used for calibration, does not show any significant difference compared to the classical extrapolation approach on the second test case, yet the behavior on a wider set of models has to be studied.

Regarding the enhancements, we can point out a faster restart of embedded solvers inside the subsystems when $C^1$ smoothness of inputs is garanteed (presented in  \ref{subsubsection:smoothness_enhancement}) at the communication times: rough discontinuities do not necessarily occur at every communication time. This can only be considered for subsystems supporting $3$\up{rd} order polynomial inputs.

\bibliographystyle{apalike}
{\small
\bibliography{Cosim}}

\end{document}